\documentclass[12pt]{article}

\usepackage{epsf,epsfig,amsmath,amsfonts}
\usepackage{graphicx}

\newcommand\Tstrut{\rule{0pt}{3ex}}            
\newcommand\TMstrut{\rule{0pt}{4ex}}         
\newcommand\TBstrut{\rule{0pt}{5ex}}          
\newcommand\Bstrut{\rule[-1.5ex]{0pt}{0pt}} 
\newcommand\BMstrut{\rule[-3ex]{0pt}{0pt}} 
\newcommand\BBstrut{\rule[-4ex]{0pt}{0pt}} 

\newcommand{\bbe}{\mathbb{E}}
\newcommand{\bbp}{\mathbb{P}}

\newcommand{\bcal}{\mathcal{B}}

\newcommand{\fcal}{\mathcal{F}}

\newcommand{\ical}{\mathcal{I}}
\newcommand{\icalo}{\mathring{\ical}}

\newcommand{\be}{\begin{equation*}}
\newcommand{\ben}{\begin{equation}}
\newcommand{\ee}{\end{equation*}}
\newcommand{\een}{\end{equation}}

\newtheorem{thm}{Theorem}

\newtheorem{rem}{Remark}

\begin{document}

\title{Necessary and sufficient conditions for the $r$-excessive
local martingales to be martingales}
\author{Mikhail Urusov\thanks{
Faculty of Mathematics, University of Duisburg-Essen,
Thea-Leymann-Strasse 9, 45127 Essen,
Germany, Email: mikhail.urusov@uni-due.de}
\ \ \ \ Mihail Zervos\thanks{Department of Mathematics,
London School of Economics, Houghton Street, London WC2A 2AE, UK,
Email: mihalis.zervos@gmail.com}}
\date{\today}
\maketitle

\begin{abstract} \noindent
We consider the decreasing and the increasing $r$-excessive
functions $\varphi_r$ and $\psi_r$ that are associated with
a one-dimensional conservative regular continuous strong
Markov process $X$ with values in an interval with endpoints
$\alpha < \beta$.
We prove that the $r$-excessive local martingale
$\bigl( e^{-r (t \wedge T_\alpha)} \varphi_r (X_{t \wedge T_\alpha})
\bigr)$ $\bigl($resp., $\bigl( e^{-r (t \wedge T_\beta)} \psi_r
(X_{t \wedge T_\beta}) \bigr) \bigr)$ is a strict local
martingale if the boundary point $\alpha$ (resp., $\beta$)
is inaccessible and entrance, and a martingale
otherwise.
\smallskip

\smallskip

\end{abstract}

\section{Introduction}

We consider a one-dimensional conservative regular
continuous strong Markov process
$X = (\Omega, \fcal, \fcal_t, \bbp_x, X_t; \, t \geq 0, \, 
x \in \ical)$ with values in an interval $\ical \subseteq
[-\infty, \infty]$ with endpoints $\alpha < \beta$ that is
open, closed or semi-open (in Remark~\ref{refs}, we
provide references for all the facts we state in what
follows).
We recall that a Markov process is called conservative if
there is no killing and a one-dimensional continuous strong
Markov process with state space $\ical$ is called regular if
\be
\bbp_x (T_y < \infty) > 0 \quad \text{for all } x \in \icalo
\text{ and } y \in \ical ,
\ee
where $\icalo = \mbox{} ]\alpha, \beta[$.
Throughout the paper we denote
\be
T_y = \inf \{ t \geq 0 \mid \ X_t = y \} , \quad \text{for }
y \in [\alpha, \beta] ,
\ee
with the usual convention that $\inf \emptyset = \infty$.
Also, we denote by $p$ and $m$ the scale function
and the speed measure of $X$.
Furthermore, we recall that the boundary point $\alpha$
(resp., $\beta$) is called inaccessible if
\ben
\bbp_x (T_\alpha < \infty) = 0 \quad \Bigl( \text{resp., }
\bbp_x (T_\beta < \infty) = 0 \Bigr) \quad \text{for all }
x \in \icalo , \label{in/accessible}
\een
and accessible otherwise.

Given any $r>0$, there exist a continuous decreasing
function $\varphi_r : \icalo \rightarrow \mbox{} ]0, \infty[$
and a continuous increasing function $\psi_r : \icalo
\rightarrow \mbox{} ]0, \infty[$ function that are determined
uniquely, up to multiplicative constants, by the expressions
\ben
\varphi_r (y) = \varphi_r (x) \, \bbe_y \left[ e^{-rT_x}
\right] \quad \text{and} \quad
\psi_r (x) = \psi_r (y) \, \bbe_x \left[ e^{-rT_y} \right]
 \quad \text{for all } x<y \text{ in } \icalo . \label{phi-psi}
\een
These functions are often called $r$-excessive.
Since they are monotone, they can be extended to
$[\alpha, \beta]$ by defining
\be
\varphi_r (\alpha) = \lim _{x \downarrow \alpha}
\varphi_r (x) , \quad 
\psi (\alpha) = \lim _{x \downarrow \alpha}
\psi_r (x) , \quad 
\varphi_r (\beta) = \lim _{x \uparrow \beta}
\varphi_r (x) \quad \text{and} \quad
\psi_r (\beta) = \lim _{x \uparrow \beta} \psi (x) .
\ee
Furthermore,
\ben
\alpha \ \Bigl(\text{resp., } \beta \Bigr) \text{ is
inaccessible if and only if } \varphi_r (\alpha) = \infty
\ \Bigl( \text{resp., } \psi_r (\beta) = \infty \Bigr) .
\label{in/accessible-exc}
\een
An important property of $\psi_r$ and $\varphi_r$ is that
\ben
\text{the processes } \Bigl( e^{-r(t \wedge T_\alpha)}
\varphi_r (X_{t \wedge T_\alpha}) \Bigr)
\text{ and } \Bigl( e^{-r (t \wedge T_\beta)} \psi_r
(X_{t \wedge T_\beta}) \Bigr) \text{ are } \bbp_x
\text{-local martingales} \label{phi-psi-lms}
\een
for all $x \in \ical$.
Despite their widespread use, these processes still do not
have a standard name.
In this paper, we refer to them as {\em $r$-excessive
$\bbp_x$-local martingales\/}.

Beyond the central role that they play in the theory of
one-dimensional diffusions, the $r$-excessive functions
$\psi_r$, $\varphi_r$ and their associated $r$-excessive
$\bbp_x$-local martingales have been used extensively in
the analysis and the solution of numerous optimal stopping
and stochastic control problems involving one-dimensional
diffusions.
This most widespread use has motivated this paper.
We refrain from trying to provide any relevant representative
references because the use of the $r$-excessive functions
and local martingales in applications
of stochastic analysis has become folklore.

We derive necessary and sufficient conditions for the
$r$-excessive $\bbp_x$-local martingales to be
$\bbp_x$-martingales.
If $\beta$ is accessible, then $\bigl( e^{-r (t \wedge T_\beta)}
\psi_r (X_{t \wedge T_\beta}) \bigr)$ is a $\bbp_x$-martingale
for all $x \in \ical$ because it is a bounded $\bbp_x$-local
martingale.
On the other hand, we prove that, if $\beta$ is inaccessible,
then ({\em i\/}) $\bigl( e^{-rt} \psi_r (X_t) \bigr)$ is a
$\bbp_x$-martingale for all $x \in \ical$ if $\beta$ is a natural
boundary point, and ({\em ii\/}) $\bigl( e^{-rt} \psi_r (X_t) \bigr)$
is a strict $\bbp_x$-local martingale for all $x \in \ical$ if
$\beta$ is an entrance boundary point, unless $\alpha$
is absorbing and $x = \alpha$, in which case the process
$\bigl( e^{-rt} \psi_r (X_t) \bigr)$ under $\bbp_\alpha$ is
identically equal to 0.
We emphasise that we do not impose any restrictions
on the boundary behaviour of $\alpha$ if this is accessible:
it can be instantaneously or slowly reflecting as well as
absorbing.
Symmetric statements hold true for the $\bbp_x$-local
martingale $\bigl( e^{-r(t \wedge T_\alpha)} \varphi_r
(X_{t \wedge T_\alpha}) \bigr)$.
We expand on these statements in Theorem~\ref{thm:main},
our main result.

A result of a closely related nature has been established
by Kotani~\cite{K}: the $\bbp_x$-local martingale $\bigl(
p(X_{t \wedge T_\alpha \wedge T_\beta}) \bigr)$
is a $\bbp_x$-martingale if and only if neither $\alpha$
nor $\beta$ is an entrance boundary point.
In fact, Delbaen and Shirakawa~\cite{DS} had earlier
established this result in a special case.
The further analysis in Hulley~\cite[Chapter~2]{H} is also
worth mentioning.
Furthermore, Gushchin, Urusov and Zervos~\cite{GUZ}
complemented this result by showing that the
$\bbp_x$-local martingale $\bigl(
p(X_{t \wedge T_\alpha \wedge T_\beta}) \bigr)$ is a
$\bbp_x$-supermaringale (resp., $\bbp_x$-submartingale)
if and only if $\alpha$ (resp., $\beta$) is not an entrance
boundary point.

\section{The main result}

Before addressing our main result, we recall that
the boundary point $\beta$ is inaccessible if and only if
\ben
\int _x^\beta m \bigl( [x, y[ \bigr) \, p(dy) = \infty ,
\label{in/accessible-an}
\een
where $x \in \icalo$ and $p(dy)$ is the atomless measure
on $(\icalo,\bcal(\icalo))$ satisfying $p\bigl( ]a,b] \bigr)
= p(b) - p(a)$ for $\alpha < a < b < \beta$
(see also the definition in (\ref{in/accessible}) as well
as (\ref{in/accessible-exc})).
We note that this characterisation does not depend on
the choice of $x \in \icalo$ because $m$ is a Radon
measure.
Also, if $\beta$ is inaccessible, then it is called
natural if
\ben
\lim _{x \downarrow \alpha} \bbp_x (T_y < t) = 0
\quad \text{for all } y \in \icalo \text{ and } t > 0
\label{natural}
\een
and entrance otherwise, namely, if
\ben
\lim _{x \downarrow \alpha} \bbp_x (T_y < t) > 0
\quad \text{for some } y \in \icalo \text{ and } t > 0 .
\label{entrance}
\een
In terms of an analytic characterisation, if $\beta$ is
inaccessible then it is
\begin{align}
\text{natural if } & \int _x^\beta m \bigl( [y, \beta[ \bigr)
\, p(dy) = \infty \label{natural-an} \\
\text{and} \quad \text{entrance if } & \int _x^\beta m
\bigl( [y, \beta[ \bigr) \, p(dy) < \infty , \label{entrance-an}
\end{align}
where the choice of $x \in \icalo$ is again arbitrary.
Furthermore, we recall that the $r$-excessive functions
$\varphi_r$ and $\psi_r$ satisfy the second order
differential equation
\be
\frac{d}{dm} \frac{d^+ f}{dp} (x) = r f(x)
\ee
in the sense that the limits
\be
\frac{d^+ f}{dp} (x) = \lim _{\varepsilon \downarrow 0}
\frac{f(x+\varepsilon) - f(x)}{p(x+\varepsilon) -
p(x)}
\ee
exist for all $x \in \icalo$ and
\ben
\frac{d^+ f}{dp} (x_2) - \frac{d^+ f}{dp} (x_1) =
r \int _{]x_1, x_2]} f(y) \, m(dy) \quad \text{for all }
\alpha < x_1 < x_2 < \beta . \label{phi-psi-2drv}
\een

\begin{rem} \label{refs} {\rm
All of the claims that we have made about the
diffusion $X$, its boundary classification and its
$r$-excessive functions are standard, and can
be found in
It\^{o} and McKean~\cite[Chapter~4]{IM},
Rogers and Williams~\cite[Section~V.7]{RW},
Karlin and Taylor~\cite[Chapter~15]{KT},
Revuz and Yor~\cite[Section~VII.3]{RY}, and
Borodin and Salminen~\cite[Chapter~II]{BS}.
In terms of boundary classification, the terminology that
we have adopted is the same as the one in Karlin and
Taylor~\cite[Table~15.6.2]{KT} and is consistent with
the one in Revuz and Yor~\cite[Section~VII.3]{RY}
and Rogers and Williams~\cite[Section~V.51]{RW}.
On the other hand, It\^{o} and McKean~\cite{IM}
use the terminology ``not exit'', ``not entrance,
not exit'' and ``entrance, not exit'' in place
of ``inaccessible'', ``natural'' and ``entrance'',
while Borodin and Salminen~\cite[Section~II.1]{BS}
use the terminology ``not exit'', ``natural'' and
``entrance-not-exit'' in place of ``inaccessible'',
``natural'' and ``entrance''.
} \mbox{} \hfill $\Box$
\end{rem}

The proof of our main result, which is captured by
(A) in the following table, involves establishing
first (B)--(D) using (E).
We state explicitly all of these cases as well as (F)
due to their independent interest as well as for
completeness.

\begin{thm} \label{thm:main}
The following statements hold true:
\smallskip

\noindent {\rm (I)}
If $\beta$ is accessible, namely, if the conditions in
(\ref{in/accessible}) and (\ref{in/accessible-an}) fail, then 
the process $\bigl( e^{-r (t \wedge T_\beta)} \psi_r
(X_{t \wedge T_\beta}) \bigr)$ is a $\bbp_x$-martingale
for all $x \in \ical$.
\smallskip

\noindent {\rm (II)}
Suppose that $\beta$ is inaccessible, namely, the conditions
in (\ref{in/accessible}) and (\ref{in/accessible-an}) hold true.
If $\beta$ is natural, namely, if the conditions in (\ref{natural})
and (\ref{natural-an}) hold true, then the process
$\bigl( e^{-rt} \psi_r (X_t) \bigr)$ is a $\bbp_x$-martingale
for all $x \in \ical$.
On the other hand, if $\beta$ is entrance, namely, if the
conditions in (\ref{entrance}) and (\ref{entrance-an}) hold
true, then the process $\bigl( e^{-rt} \psi_r (X_t) \bigr)$ is a
strict $\bbp_x$-local martingale for all $x \in \ical$, unless
$\alpha$ is absorbing and $x = \alpha$.
Furthermore, the equivalences suggested by the following
table hold true (note that all limits appearing here indeed exist).

\begin{rm}
\begin{center}
\begin{tabular}{|c||c|c|}
\hline
& $\beta$ is natural & $\beta$ is entrance
\Tstrut\Bstrut \\ \hline\hline
(A) & $\forall r>0, \, \forall x \in \ical , \ \bigl( e^{-rt} \psi_r
(X_t) \bigr)$ &
$\forall r>0, \, \forall x \in \icalo , \ \bigl( e^{-rt} \psi_r
(X_t) \bigr)$
\Tstrut \\
& is a $\bbp_x$-martingale & is a strict $\bbp_x$-local martingale
\Bstrut \\ \hline
(B) & $\displaystyle \forall s>r>0, \ \lim _{x \uparrow \beta}
\frac{\psi_s (x)}{\psi_r (x)} = \infty$
& $\displaystyle \forall s>r>0, \ \lim _{x \uparrow \beta}
\frac{\psi_s (x)}{\psi_r (x)} \in \mbox{} ]0, \infty[$
\TMstrut\BMstrut \\ \hline
(C) & $\displaystyle \forall r>0, \ \lim _{x \uparrow \beta}
\frac{\psi_r (x)}{p(x)} = \infty$
& $\displaystyle \forall r>0, \ \lim _{x \uparrow \beta}
\frac{\psi_r (x)}{p(x)}  \in \mbox{} ]0, \infty[$
\TMstrut\BMstrut \\ \hline
(D) & $\displaystyle \forall s>r>0, \ \lim _{x \uparrow \beta}
\frac{\frac{d^+ \psi _s}{dp} (x)}{\frac{d^+ \psi _r}{dp} (x)} = \infty$
& $\displaystyle \forall s>r>0, \ \lim _{x \uparrow \beta}
\frac{\frac{d^+ \psi _s}{dp} (x)}{\frac{d^+ \psi _r}{dp} (x)} \in
\mbox{} ]0, \infty[$
\TBstrut\BBstrut \\ \hline 
(E) & $\displaystyle \forall r>0, \ \lim _{x \uparrow \beta}
\frac{d^+ \psi _r}{dp} (x) = \infty$& $\displaystyle \forall r>0,
\ \lim _{x \uparrow \beta} \frac{d^+ \psi _r}{dp} (x)
\in \mbox{} ]0, \infty[$
\TMstrut\BMstrut \\ \hline 
(F) & $\displaystyle \forall r>0, \, \forall x \in \icalo , \
\int _{[x,\beta[} \psi_r (y) \, m(dy) = \infty$ &
$\displaystyle \forall r>0, \, \forall x \in \icalo , \
\int _{[x,\beta[} \psi_r (y) \, m(dy) < \infty$
\TMstrut\BMstrut \\ \hline
\end{tabular}
\label{Table2}
\end{center}
\end{rm}
\smallskip

\noindent {\rm (III)}
Symmetric results hold true for the process
$\bigl( e^{-r (t \wedge T_\alpha)} \varphi_r
(X_{t \wedge T_\alpha}) \bigr)$.
\end{thm}
{\bf Proof.}
Statement (I) follows immediately because
$\bigl( e^{-r (t \wedge T_\beta)} \psi_r (X_{t \wedge T_\beta})
\bigr)$ is a bounded $\bbp_x$-local martingale (see also
(\ref{in/accessible-exc})).
To prove~(II), we assume in what follows that $\beta$ is
inaccessible, which implies that
\ben
\lim _{x \rightarrow \beta} \psi_r (x) = \infty \quad \text{for all }
r>0 . \label{lim-psi}
\een
The results in (E) and (F) appear in the fourth and the sixth
row of Table~1 in It\^{o} and McKean~\cite[Section~4.6]{IM}
(see the third and fourth columns of that table; also, note
that (F) follows immediately from (E) and (\ref{phi-psi-2drv})).
Also, (C) follows from (E) and the calculation
\begin{align}
\lim _{x \uparrow \beta} \frac{\psi_r (x)}{p(x)} & =
\lim _{y \uparrow p(\beta)} \frac{\psi_r \bigl( p^{-1} (y)
\bigr)}{y} = \lim _{y \uparrow p(\beta)} \frac{d^+ \psi_r
\circ p^{-1}}{dy} (y) = \lim _{x \uparrow \beta}
\frac{d^+ \psi_r}{dp} (x) , \nonumber
\end{align}
in which we have used L'H\^{o}pital's rule.

We now show that
\ben
\text{the limits } \lim _{x \uparrow \beta} \frac{\psi_r (x)}
{\psi_s (x)} \text{ exist in } [0, \infty[ \text{ for all }
s > r > 0 \label{lim-psi-rs}
\een
as well as that (A) and (B) are equivalent.
To this end, we consider
an initial condition $x \in \mathring{\ical}$, a point
$\bar{\beta} \in \mbox{} ]x, \beta[$ and constant $s>r>0$,
and we use the integration by parts formula to calculate
\be
e^{-(s-r)t} M_{t \wedge T_{\bar{\beta}}} = \psi_r (x) -
(s-r) \int _0^{t \wedge T_{\bar{\beta}}} e^{-(s-r)u} M_u
\, du + \int _0^{t \wedge T_{\bar{\beta}}} e^{-(s-r)u}
\, dM_u ,
\ee
where $M_t = e^{-rt} \psi_r (X_t)$.
The process $(M_{t \wedge T_{\bar{\beta}}}, \ t \geq 0)$ is
a $\bbp_x$-square integrable martingale because it
is a bounded $\bbp_x$-local martingale.
Therefore, the stochastic integral in this identity has zero
expectation.
In view of this observation and the dominated and
monotone convergence theorems, we can see that
\begin{align}
\psi_r (\bar{\beta}) \, \bbe_x \left[ e^{-s T_{\bar{\beta}}} \right]
& = \lim _{t \rightarrow \infty} \bbe_x \left[ e^{-s (t \wedge
T_{\bar{\beta}})} \psi_r (X_{t \wedge T_{\bar{\beta}}})
\right] \nonumber \\
& = \lim _{t \rightarrow \infty} \bbe_x \left[ e^{-(s-r)
(t \wedge T_{\bar{\beta}})} M_{t \wedge T_{\bar{\beta}}}
\right] \nonumber \\
& = \psi_r (x) - (s-r) \, \bbe_x \left[ \int _0^{T_{\bar{\beta}}}
e^{-(s-r)u} M_u \, du \right] \nonumber \\
& = \psi_r (x) - (s-r) \, \bbe_x \left[ \int _0^{T_{\bar{\beta}}}
e^{-su} \psi_r (X_u) \, du \right] . \nonumber
\end{align}
Combining this calculation with the definition of $\psi_s$
as in (\ref{phi-psi}), we obtain
\be
\psi_r (\bar{\beta}) \frac{\psi_s (x)}{\psi_s (\bar{\beta})}
= \psi_r (x) - (s-r) \, \bbe_x \left[ \int _0^{T_{\bar{\beta}}}
e^{-su} \psi_r (X_u) \, du \right] .
\ee
In view of the monotone convergence theorem and the
assumption that $\beta$ is inaccessible, it follows that
\ben
\lim _{\bar{\beta} \uparrow \beta} \frac{\psi_r (\bar{\beta})}
{\psi_s (\bar{\beta})} = \frac{\psi_r (x)}{\psi_s (x)} -
\frac{s-r}{\psi_s (x)} \int _0^\infty e^{-(s-r)u} \, \bbe_x
\bigl[ e^{-ru} \psi_r (X_u) \bigr] \, du . \label{psir/psis}
\een
This identity and the positivity of $\psi_r$ imply
that (\ref{lim-psi-rs}) is indeed true.
Furthermore, since the process $\bigl( e^{-rt} \psi_r
(X_t) , \ t \geq 0 \bigr)$ a positive $\bbp_x$-local
martingale, it is a $\bbp_x$-supermartingale.
Therefore,
\be
\bbe_x \bigl[ e^{-rt} \psi_r (X_t) \bigr] \leq 
\psi_r (x) \quad \text{for all } t \geq 0 ,
\ee
with equality holding if and only if $\bigl(e^{-rt} \psi_r
(X_t) , \ t \geq 0 \bigr)$ is a $\bbp_x$-martingale.
In view of this observation, we can see that
(\ref{psir/psis}) implies that
$\lim _{\bar{\beta} \uparrow \beta} \psi_r (\bar{\beta})
/ \psi_s (\bar{\beta}) = 0$ if and only if $(e^{-rt} \psi_r
(X_t) , \ t \geq 0)$ is a $\bbp_x$-martingale,
and the equivalence of (A) and (B) follows.

To complete the proof, we need to establish (B) and (D).
To this end, we note that (\ref{lim-psi}), (\ref{lim-psi-rs})
and L'H\^{o}pital's rule imply that
\be
\lim _{x \uparrow \beta} \frac{\psi_r (x)}{\psi_s (x)}
= \lim _{y \uparrow p(\beta)} \frac{\psi_r \bigl(
p^{-1} (y) \bigr)}{\psi_s \bigl( p^{-1} (y) \bigr)}
= \lim _{y \uparrow p(\beta)} \frac{\frac{d^+ \psi_r
\circ p^{-1}}{dy} (y)}{\frac{d^+ \psi_s \circ
p^{-1}}{dy} (y)} = \lim _{x \uparrow \beta}
\frac{\frac{d^+ \psi_r}{dp} (x)}{\frac{d^+ \psi_s}{dp} (x)}
\ee
whenever the last limit exists.
If $\beta$ is an entrance boundary point, then
this calculation and the corresponding statement in
(E) imply that the corresponding claims in (B) and
(D) are indeed true.

On the other hand, if $\beta$ is a natural boundary
point, then we can use (\ref{phi-psi-2drv}), (\ref{lim-psi})
and (\ref{lim-psi-rs}) to see that, given any $x_1 \in \icalo$,
\begin{align}
\lim _{x \uparrow \beta} \frac{\psi_r (x)}{\psi_s (x)}
& = \lim _{x \uparrow \beta} \frac{\frac{d^+ \psi_r}{dp} (x)}
{\frac{d^+ \psi_s}{dp} (x)}
= \lim _{x \uparrow \beta} \frac{\frac{d^+ \psi_r}{dp}
(x_1) + r \int _{]x_1, x]} \psi_r (y) \, m(dy)}
{\frac{d^+ \psi_s}{dp} (x_1) + s \int _{]x_1, x]} \psi_s (y)
\, m(dy)} \nonumber \\
& = \frac{r}{s} \lim _{x \uparrow \beta} \frac{
\int _{]x_1, x]} \psi_r (y) \, m(dy)} {\int _{]x_1, x]}
\psi_s (y) \, m(dy)}
= \frac{r}{s} \lim _{x \uparrow \beta} \frac{\psi_r (x)}
{\psi_s (x)} . \nonumber
\end{align}
In view of (\ref{lim-psi-rs}), all these limits are equal to 0,
and the corresponding claims in (B) and (D) follow.
\mbox{}\hfill$\Box$

\end{document}